# GOLDBACH CONJECTURE AND FIRST-ORDER ARITHMETIC

## BY FERNANDO REVILLA

Department of Mathematics. I.E.S. Sta Teresa de Jesús.
Comunidad de Madrid. Spain.
frej0002@ficus.pntic.mec.es

ABSTRACT. Using the concepts of Hyperbolic Classification of Natural Numbers, Essential Regions and Goldbach Conjecture Function ( [1] ) we prove that the existence of a proof of the statement "*Every even integer greater than 2 can be expressed as a sum of two primes*" in First-Order Arithmetic ( $\mathcal{N}$ ), would imply the existence of another proof in a certain extension $\mathcal{N}^*$ of $\mathcal{N}$ that would not be valid in all states of time associated to natural numbers by means of adequate dynamic processes.

**1. First-Order Arithmetic** Consider the first order language $\mathcal{L}$ such that its alphabet of symbols is: ($\mathcal{L}$1) $x_1, x_2, ...$ (variables) ($\mathcal{L}$2) $a_1, a_2, ...$ (individual constants) ($\mathcal{L}$3) $A_1^1, A_2^1, ..., A_1^2, A_2^2, ..., A_1^3, A_2^3, ...$ (predicative letters) ($\mathcal{L}$4) $f_1^1, f_2^1, ..., f_1^2, f_2^2, ..., f_1^3, f_2^3,$ ... (function letters)  ($\mathcal{L}$5) $(,),,$ (punctuation signs)  ($\mathcal{L}$6) $\neg, \rightarrow$ (negation and conditional connectives) ($\mathcal{L}$7) $\forall$ (universal quantifier)

Occasionally we will use the symbols $\exists$, $\wedge$, $\vee$ and $\leftrightarrow$ :

$$(\exists x_i)\mathcal{A} \quad \text{is an abbreviation of} \quad \left(\neg\left((\forall x_i)(\neg\mathcal{A})\right)\right)$$

$$(\mathcal{A} \wedge \mathcal{B}) \quad \text{is an abbreviation of} \quad \left(\neg(\mathcal{A} \rightarrow (\neg\mathcal{B}))\right)$$

$$(\mathcal{A} \vee \mathcal{B}) \quad \text{is an abbreviation of} \quad \left((\neg\mathcal{A}) \rightarrow \mathcal{B}\right)$$

$$\mathcal{A} \leftrightarrow \mathcal{B} \quad \text{is an abbreviation of} \quad (\mathcal{A} \rightarrow \mathcal{B}) \wedge (\mathcal{B} \rightarrow \mathcal{A})$$

Let $\mathcal{L}$ be a first order language and a formal deductive system $K_\mathcal{L}$ defined by the following axioms and rules of deduction:

*Axioms.*

Let $\mathcal{A}, \mathcal{B}, \mathcal{C}$ be any well formed formulas (*wff*) of $\mathcal{L}$. The following are axioms of $K_\mathcal{L}$
(K1) $(\mathcal{A} \rightarrow (\mathcal{B} \rightarrow \mathcal{A}))$   (K2) $(\mathcal{A} \rightarrow (\mathcal{B} \rightarrow \mathcal{C})) \rightarrow ((\mathcal{A} \rightarrow \mathcal{B}) \rightarrow (\mathcal{A} \rightarrow \mathcal{C}))$.
(K3)$(\neg\mathcal{A} \rightarrow \neg\mathcal{B}) \rightarrow (\mathcal{B} \rightarrow \mathcal{A})$  (K4)$((\forall x_i)\mathcal{A} \rightarrow \mathcal{A})$, if $x_i$ does not occur free in $\mathcal{A}$.
(K5)$((\forall x_i)\mathcal{A}(x_i) \rightarrow \mathcal{A}(t))$, if $\mathcal{A}(x_i)$ is a *wff* of $\mathcal{L}$ and $t$ is a term in $\mathcal{L}$ which is free for $x_i$ in $\mathcal{A}(x_i)$   (K6) $(\forall x_i)(\mathcal{A} \rightarrow \mathcal{B}) \rightarrow (\mathcal{A} \rightarrow (\forall x_i)\mathcal{B})$, if $\mathcal{A}$ contains no free occurrence of the variable $x_i$.

*Rules*

(1) *Modus Ponens,* i.e. from $\mathcal{A}$ and $\mathcal{A} \to \mathcal{B}$ deduce $\mathcal{B}$, where $\mathcal{A}$ and $\mathcal{B}$ are any *wff* of $\mathcal{L}$.

(2) *Generalization*, i.e. from $\mathcal{A}$ deduce $(\forall x_i)\mathcal{A}$, where $\mathcal{A}$ is any *wff* of $\mathcal{L}$ and $x_i$ is any variable.

The language $\mathcal{L}_\mathbb{N}$ we take to include the variables $x_1, x_2, \ldots$, the individual constant $a_1$ (for 0), the function letters $f_1^1, f_1^2, f_2^2$ (successor, sum and product) and the predicative symbol =, as well as punctuation, connectives and universal quantifier. Let us denote by $\mathcal{N}$ the first order system which is the extension of $K_{\mathcal{L}_\mathbb{N}}$ obtained by including the following six axioms ( (N1), …, (N6) ) and one axiom scheme ( (N7) ).

(N1) $(\forall x_1) \neg (f_1^1(x_1) = a_1)$. (N2) $(\forall x_1)(\forall x_2)(f_1^1(x_1) = f_1^1(x_2) \to x_1 = x_2)$.

(N3) $(\forall x_1)(f_1^2(x_1, a_1) = x_1)$. (N4) $(\forall x_1)(\forall x_2)(f_1^2(x_1, f_1^1(x_2)) = f_1^1(f_1^2(x_1, x_2)))$.

(N5) $(\forall x_1)(f_2^2(x_1, a_1) = a_1)$. (N6) $(\forall x_1)(\forall x_2)(f_2^2(x_1, f_1^1(x_2)) = f_1^2(f_2^2(x_1, x_2), x_1))$.

(N7) $\mathcal{A}(a_1) \to \left( (\forall x_1)(\mathcal{A}(x_1) \to \mathcal{A}(f_1^1(x_1))) \to (\forall x_1)\mathcal{A}(x_1) \right)$ (for each *wff* $\mathcal{A}(x_1)$ of $\mathcal{L}_\mathbb{N}$ in which $x_1$ occurs free). (First-Order Arithmetic).

**2. Family of interpretations of $\mathcal{N}$** Consider the set:
$$\mathfrak{N} = \{ \alpha \in \mathbb{N} : (\alpha \text{ even}) \wedge (\alpha \geq 16) \wedge (\frac{\alpha}{2} \text{ non-prime}) \wedge (\alpha - 3 \text{ non-prime}) \}.$$

Consider the family of $\mathbb{R}^+$ coding functions:
$$\{ \psi_{(\alpha, u)} : \mathbb{R}^+ \to [0, M_{\psi_{(\alpha, u)}}) : (\alpha, u) \in \mathfrak{N} \times (1, +\infty) \}$$

such that $\xi_0 = \xi_1 = \xi_{\alpha/2} = \xi_i = 1$ ($\forall i \geq \alpha - 4$) and the remaining coefficients have been chosen in such a way as to create the corresponding scalar Goldbach Conjecture function $\mathfrak{G} = \mathfrak{G}_{(\alpha, 1, 1, u)}$ associated to $\alpha$. ( [1] (3.1.13) ). For every natural number *n*, we denote $\hat{n}_{(\alpha, u)} := \psi_{(\alpha, u)}(n)$ and we define the set: $\widehat{\mathbb{N}}_{(\alpha, u)} := \{ \hat{n}_{(\alpha, u)} : n \in \mathbb{N} \}$. We construct the family $\{ I_{(\alpha, u)} : (\alpha, u) \in \mathfrak{N} \times (1, +\infty) \}$ of interpretations of $\mathcal{N}$ in the following way:

(a) *Domain of $I_{(\alpha, u)}$*: $D_{I_{(\alpha, u)}} = \widehat{\mathbb{N}}_{(\alpha, u)}$.

(b) *Distinguished elements*: $(\overline{a_1})_{(\alpha, u)} = 0$, $(\overline{a_2})_{(\alpha, u)} = \hat{1}_{(\alpha, u)}$. (c) *Functions on $\widehat{\mathbb{N}}_{(\alpha, u)}$*:

($c_1$) Successor function on $\widehat{\mathbb{N}}_{(\alpha, u)}$:
$$\left(\overline{f_1^1}\right)_{(\alpha, u)} : \widehat{\mathbb{N}}_{(\alpha, u)} \to \widehat{\mathbb{N}}_{(\alpha, u)}, \quad \left(\overline{f_1^1}\right)_{(\alpha, u)}(x_1) = x_1 \oplus_{(\alpha, u)} \hat{1}_{(\alpha, u)}$$

( where $\oplus_{(\alpha, u)}$ represents the usual sum of $\mathbb{N}$ transported to $\widehat{\mathbb{N}}_{(\alpha, u)}$ through $\psi_{(\alpha, u)}$ )

($c_2$) Sum function on $\widehat{\mathbb{N}}_{(\alpha, u)}$:
$$\left(\overline{f_1^2}\right)_{(\alpha, u)} : \left(\widehat{\mathbb{N}}_{(\alpha, u)}\right)^2 \to \widehat{\mathbb{N}}_{(\alpha, u)}, \quad \left(\overline{f_1^2}\right)_{(\alpha, u)}(x_1, x_2) = x_1 \oplus_{(\alpha, u)} x_2$$

($c_3$) Product function on $\widehat{\mathbb{N}}_{(\alpha,u)}$:

$$\left(\overline{f_2^2}\right)_{(\alpha,u)} : \left(\widehat{\mathbb{N}}_{(\alpha,u)}\right)^2 \to \widehat{\mathbb{N}}_{(\alpha,u)}, \quad \left(\overline{f_2^2}\right)_{(\alpha,u)}(x_1, x_2) = x_1 \otimes_{(\alpha,u)} x_2$$

(where $\otimes_{(\alpha,u)}$ represents the usual product of $\mathbb{N}$ transported to $\widehat{\mathbb{N}}_{(\alpha,u)}$ through $\psi_{(\alpha,u)}$)

**3. Particular interpretation of $\mathcal{N}$** We add to the aforementioned interpretations, a particular one, in the following manner. For every $\alpha \in \mathfrak{N}$ we verify:

$$\lim_{u\to 1^+} \widehat{n}_{(\alpha,u)} = \lim_{u\to 1^+} \psi_{(\alpha,u)}(n) = \psi_{(\alpha,1)}(n) = n$$

So, we obtain the following interpretation of $\mathcal{N}$:

(a)' *Domain of* $I_{(\alpha,1)}$: $D_{I_{(\alpha,1)}} = \widehat{\mathbb{N}}_{(\alpha,1)} = \mathbb{N}$, because as a consequence of [1] (3.2.2) all $\xi_i$ coefficients are equal to 1 and $\psi_{(\alpha,1)} = id : \mathbb{R}^+ \to \mathbb{R}^+$ (*id*: identity function)

(b)' *Distinguished elements*: $\left(\overline{a_1}\right)_{(\alpha,u)} = 0$, $\left(\overline{a_2}\right)_{(\alpha,1)} = \widehat{1}_{(\alpha,1)} = 1$.

(c)' *Functions on* $\widehat{\mathbb{N}}_{(\alpha,1)}$:

($c_1$)' Successor function on $\mathbb{N}$:

$$\left(\overline{f_1^1}\right)_{(\alpha,1)} : \mathbb{N} \to \mathbb{N}, \quad \left(\overline{f_1^1}\right)_{(\alpha,1)}(x_1) = x_1 + 1$$

(where + represents the usual sum of $\mathbb{N}$)

($c_2$)' Sum function on $\widehat{\mathbb{N}}_{(\alpha,1)}$:

$$\left(\overline{f_1^2}\right)_{(\alpha,1)} : \mathbb{N}^2 \to \mathbb{N}, \quad \left(\overline{f_1^2}\right)_{(\alpha,1)}(x_1, x_2) = x_1 + x_2$$

($c_3$)' Product function on $\widehat{\mathbb{N}}_{(\alpha,1)}$:

$$\left(\overline{f_2^2}\right)_{(\alpha,1)} : \mathbb{N}^2 \to \mathbb{N}, \quad \left(\overline{f_2^2}\right)_{(\alpha,1)}(x_1, x_2) = x_1 \cdot x_2$$

(where $\cdot$ represents the usual product of $\mathbb{N}$)

**4. Consequence of the existence of a proof of the Goldbach Conjecture in First-Order Arithmetic.** Denote by $\mathcal{G}$ the *wff* in $\mathcal{N}$ which corresponds to the sentence:

"For every $\alpha$ even number such that $\alpha \geq 16$ and $\alpha - 3$ non-prime and $\dfrac{\alpha}{2}$ non-prime, $\alpha$ can be expressed as a sum of two primes". In a certain extension $\mathcal{N}^*$ of $\mathcal{N}$ denote by $\mathcal{E}_u$ ($u > 1$) the *wff* which corresponds to the sentence: " For every $\alpha$ even number such that $\alpha \geq 16$ and $\alpha - 3$ non-prime and $\dfrac{\alpha}{2}$ non-prime, two essential points associated to $\alpha$ by means of $\psi_{(\alpha,u)}$ are repeated".

Suppose there exists a proof of the Goldbach Conjecture in the formal system $\mathcal{N}$, then there exists in a sequence of *wff* $\mathcal{A}_1$, ..., $\mathcal{G}$ of $\mathcal{N}$ such that each member is an axiom of $\mathcal{N}$ or follows from previous members of the sequence by Modus Ponens or Generalization. According to [1] (2.8.3 and 3.2.2), in the extension $\mathcal{N}^*$:
(i) For every $u > 1$, $\left(\mathcal{G} \leftrightarrow \mathcal{E}_u\right)$ is a theorem of $\mathcal{N}^*$.
(ii) For $u = 1$ all the essential points are:

$$P_i^{(\alpha,1)} = \left(\frac{1}{2}, -\frac{1}{2}\right) \ (\forall i \in \mathbb{N}: 4 \leq i \leq \frac{\alpha}{2} - 1)$$

Then, through (i) above, for every $u$ belonging to $(1, +\infty)$ the sequence

$$\mathcal{A}_1, \ldots, \mathcal{G}, \mathcal{E}_u$$

is a proof of the Goldbach Conjecture in the formal system $\mathcal{N}^*$. We now take limits where $u \to 1^+$ that is, we change the states of time associated to natural numbers by means of the dynamic process created by $\psi_{(\alpha,u)}$. Then, we obtain the sequence:

$$\mathcal{A}_1, \ldots, \mathcal{G}, \mathcal{E}_1$$

But as a consequence of (ii) above, this sequence is not a proof of the mentioned conjecture. That is, the existence of a proof of the Goldbach Conjecture in First-Order Arithmetic would imply the existence of another proof in a certain extension of $\mathcal{N}$ that would not be valid in all states of time associated to natural numbers created by means of adequate dynamic processes.